\documentclass[12pt,a4paper]{article}

\usepackage{amsmath}
\usepackage{amssymb}
\usepackage{latexsym}
\usepackage{srcltx}
\usepackage{color}
\usepackage{graphicx}

\newcommand{\re}{{\mathbb R}}
\newcommand{\n}{{\mathbb N}}

\newcommand{\cA}{{\mathcal{A}}}

\newcommand{\cQ}{{\mathcal{Q}}}

\newcommand{\cM}{{\mathcal{M}}}

\newcommand{\by}{{\boldsymbol y}}
\newcommand{\bx}{{\boldsymbol x}}
\newcommand{\bt}{{\boldsymbol t}}
\newcommand{\ba}{{\boldsymbol a}}
\newcommand{\bb}{{\boldsymbol b}}
\newcommand{\bc}{{\boldsymbol c}}

\newcommand{\bz}{{\boldsymbol z}}

\newcommand{\bk}{{\boldsymbol k}}
\newcommand{\bl}{{\boldsymbol l}}
\newcommand{\bm}{{\boldsymbol m}}
\newcommand{\bn}{{\boldsymbol n}}
\newcommand{\bss}{{\boldsymbol s}}

\newtheorem{theorem}{Theorem}
\newtheorem{prop}{Proposition}
\newtheorem{lemma}{Lemma}
\newtheorem{cor}{Corollary}

\newtheorem{defi}{Definition}

\date{}

\author{E.\,A.~Avksentyev and V.\,Yu.~Protasov
\thanks{Dept. of Mechanics and Mathematics of Moscow State University and Faculty of Computer
Science of National Research University Higher School of Economics, {e-mail: \tt\small avksentjev@mail.ru, v-protassov@yandex.ru}}}
\title{Universal measure for Poncelet-type theorems
\thanks{The second author is supported by the RFBR grants Nos. 14-01-00332 and 16-04-00832, and by the grant of Dynasty foundation}}

\begin{document}
\maketitle

\begin{abstract}
We give a simple proof of the Emch closing theorem by introducing
a new invariant measure on the circle. Special cases of that measures are well-known and have been used in the
literature to prove  Poncelet's and Zigzag theorems. Some further generalizations are also obtained by applying
the new measure.

\smallskip

\smallskip

\begin{flushright}
\noindent  \textbf{AMS 2010} {\em subject
classification:  53A04, 28A25, 51N15}
\end{flushright}

\end{abstract}
\bigskip

\begin{center}
\textbf{1. Introduction}
\end{center}
\medskip

Invariant measures on circles and conics provide powerful tools in the study of
closing theorems such as Poncelet's and Steiner's porisms, Zigzag theorem, etc. We present an elementary
formula for a universal measure that generalizes several well-known invariant measures.
We show that every pair of circles generates a function on the plane, whose restriction to an arbitrary circle
 defines an invariant measure on it. Remarkable properties of that measure give new results as well as new proofs of known facts.

The {\em Poncelet closing theorem} discovered in 1813 and published in 1822~\cite{Pons}
 states that if for two circles (or quadrics)
$\alpha$ and $\delta$,  there is an $n$-sided polygon $\bx_1\ldots \bx_n$
inscribed in~$\delta$ and circumscribed
around~$\alpha$ (i.e., the straight lines containing its sides are tangent to~$\alpha$),
then there exist infinitely many such polygons, and its  vertex $\bx_1$ can be chosen on~$\delta$ arbitrarily,
provided $\bx_1 \notin \alpha$.

There are several methods to prove Poncelet's theorem. They are all nontrivial and based on various
ideas~\cite{BB, B, DR, F, HH}.
The invariant measure approach originated with Jacobi in 1828, then improved by Bertrand, and developed further in~\cite{A, K, Sc}, etc., gives an elegant and natural proof. Consider, for example, the case when the circle $\alpha$ lies inside $\delta$.
Suppose there is a measure $m(\cdot)$ on~$\delta$ such that all oriented  arcs $\stackrel{\smile }{\bx\by} \, \subset \, \delta$
whose chords touch the circle~$\alpha$ have the same value $m(\stackrel{\smile }{\bx\by}) = \tilde m$. Then the Poncelet $n$-gon exists
if and only if the number  $n\, \tilde m$ is an integer multiple of $m(\delta)$.  Since this property
does not depend on the location of the first vertex of the polygon, the Poncelet theorem for two circles follows.

For arbitrary circles $\alpha$ and $\delta$, a measure on~$\delta$ is called invariant, if
its density $\rho = m'$ satisfies the equality $\rho (\bx) |d\bx| = \rho(\by)|d\by| $, where $d\bx , d\by$ are oriented lengths
of small arcs after perturbation of an arbitrary chord~$\bx \by$ touching~$\alpha$. If a function $\rho: \delta \to \re_+$ possess this property, then $m(\stackrel{\smile }{\bx\by}) = \int_{\bx}^{\by} \rho(\bss) d\bss$
is an invariant measure (the integration is over the arc~$\stackrel{\smile }{\bx\by}$).
For arbitrary circles $\alpha$ and $\delta$, such a measure is readily available
by the formula $\rho(\bx) = 1/\sqrt{|f(\bx)|}$, where $f(\bx) = |\bx - \bc|^2 - r^2$
is the power with respect to the circle~$\alpha$ of radius $r$ centered at $\bc \in \re^2$.
This is the {\em Jacobi-Bertrand measure}. Moreover, as it was observed by Khovansky  (see~\cite{A} for an overview), if we consider an arbitrary quadratic
polynomial $f(\bx), \, \bx \in \re^2$, then the same formula also defines an invariant measure,
which  proves Poncelet's theorem for the circle $\delta$ and the quadric~$\alpha =
 \{\bx \in \re^2 \ | f(\bx) = 0\}$.  By a suitable projective transform, this leads to
 the general case of two quadrics.

In 1974 Black, Howland, and Howland~\cite{BHH} found an invariant measure for another
well-known closing theorem:

{\em Zigzag theorem.} If for given circles $\alpha, \delta$ and for a number $l > 0$,
there is a polygon with $2n$ sides all of length~$l$, with odd vertices (i.e., vertices $\bx_k$ with odd~$k$) on $\delta$ and even vertices on $\alpha$, then there exist infinitely many such polygons, and its  vertex $\bx_1$ can be chosen on~$\delta$ arbitrarily,
provided the distance from $\bx_1$ to $\alpha$ is smaller than~$l$.

Thus, if a grasshopper jumps from one circle to the other making a closed walk after $2n$ jumps,
then his walk from any point of the first circle closes after $2n$ steps, provided he can make the first jump.
This theorem was established by Emch in 1901~\cite{B}, then rediscovered by Bottema in 1965 \cite{B},
and in 1974 in~\cite{BHH}. It holds for two circles in the space as well, but we consider only the plane version.

We mention also the third popular closing theorem, the {\em Steiner theorem.} Given two circles $\alpha_0, \alpha_1$, one inside the other. Circles $\{\omega_k\}_{k \in \n}$ inscribed in the
annulus between $\alpha_0$ and $\alpha_1$ touches each other in succession ($\omega_{k}$ and $\omega_{k+2}$
are different and both tangent to $\omega_{k+1}\, , \ k \in \n$). If this series closes after $n$ steps,
i.e., $\omega_{n+1} = \omega_1$, then it does for an arbitrary initial  circle~$\omega_1$.

Those three closing   theorems are actually special cases of the Emch theorem on circular series~\cite{E}.
To formulate it we need to introduce some notation. The tangency
of two circles is called interior if one of the circles lies
inside the other. Suppose
 $\alpha_0, \alpha_1$ are circles on the plane; then for an arbitrary
 circle $\omega$ touching both $\alpha_0$ and $\alpha_1$ the {\em
index of tangency} is $0$ if there is an even number of interior
tangencies among the two ones: $\omega$ with $\alpha_0\, $ and $\,
\omega$ with $\alpha_1$. If this number is odd, then the index
is~$1$. For $i = 0,1$, let $\cM_i$ denotes the family of circles touching $\alpha_0,
\alpha_1$ with index~$i$.

Let $\alpha_0, \alpha_1$ and $\delta$ be an arbitrary triple of circles on the plane.
Choose some $i \in \{0,1\}$ and take the family $\cM_i$ of circles touching $\alpha_0, \alpha_1$
with index~$i$.
We assume $\delta \notin \cM_i$. Take arbitrary circle $\omega_1 \in \cM_i$ that intersects $\delta$
at two points $\bx_1, \bx_2$. The family $\cM_i$ contains two circles passing through $\bx_2$;
One of them is $\omega_1$, take the other and denote it $\omega_2$. The circles $\omega_2$ and $\delta$
have two points of intersection; one of them is $\bx_2$, take the other and denote it $\bx_3$, etc.
This way we obtain a {\em circular series} $\{\omega_k\}_{k \in \n}$. Each $\omega_k$ touches $\alpha_0$ and $\alpha_1$
with index~$i$ and meets the circle $\delta$ at points $\bx_k$ and $\bx_{k+1}$. This series closes after $n$ steps if
$\omega_{n+1} = \omega_1$.
 \smallskip

\begin{figure}[htb]
\center
\includegraphics[width=0.4\textwidth]{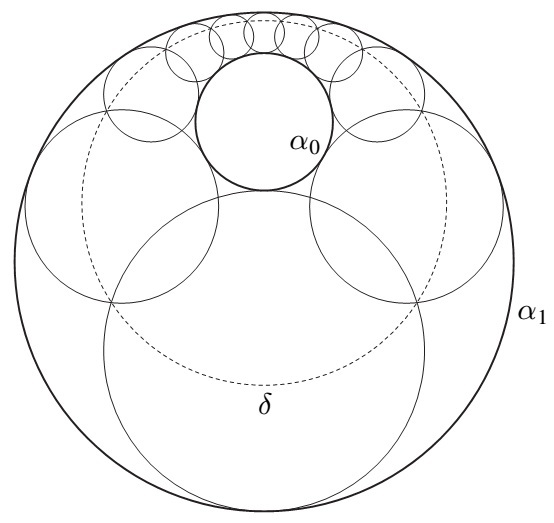}
\caption{{\footnotesize Theorem of Emch}}
\label{f1}       
\end{figure}

 \noindent \textbf{Theorem of Emch}~\cite{E}.  {\em Let $\alpha_0, \alpha_1$ and $\delta$ be arbitrary circles,
  $i \in \{0,1\}$, and $\delta \notin \cM_i$. If for some initial circle $\omega_1 \in \cM_i$,
  the circular series closes after $n$
  steps, then it does for arbitrary $\omega_1 \in \cM_i$.}
 \smallskip

This theorem was formulated by A.Emch in 1901~\cite{E}, but he gave a proof only for
non-intersecting circles $\alpha_0$ and $\alpha_1$. Probably, the author missed this aspect.
It was not before 1996 that the proof for general position of circles $\alpha_0, \alpha_1, \delta$ was given
 in~\cite{BB}. In~\cite{P1} the Emch theorem was deduced from Poncelet's theorem for quadrics, in~\cite{P2}
 an elementary geometrical proof was found, see also~\cite{A2}.

All the three famous  closing theorems follow directly from the Emch theorem.
If the radius of $\alpha_1$ tends to infinity, then we obtain in the limit the Poncelet theorem for circles
$\alpha_0$ and $\delta$. If $\alpha_0$ and $\alpha_1$ are concentric circles, then we obtain the Zigzag theorem.
Finally, if~$\delta$ is a locus of points of tangency for pairs of circles from~$\cM_i$, then we come to the Steiner theorem. See~\cite{P1} for more details.

A natural question arises, if the Emch theorem admits a proof by an invariant measure? We show that such a measure exists and, moreover, is explicitly given by a simple formula.
For an arbitrary pair of circles $\alpha_0$ and $\alpha_1$,  we consider the function
$\rho(\bx) = \frac{1}{\sqrt{|f_0(\bx)f_1(\bx)|}}$ on the plane~$\re^2$, where $f_j$ is a power
w.r.t. the circle $\alpha_j, \, j = 0,1$. In Theorem~\ref{th10} we show that this function defines
an invariant measure on any circle~$\delta \, \subset \, \re^2$. This gives a geometric proof for
the Emch theorem. Both the Jacobi-Bertrand measure and the Black-Howland
measure are special cases of this measure $\rho(\cdot)$. Therefore, it can be considered as a universal
measure for Poncelet type theorems. Simple algebraic manipulations with the formula for $\rho(\bx)$
give generalizations of Emch's theorem to pencils of circles (Section~4), to a cyclic instead of two circles
   (Section~5) and prove the equivalence of Emch's theorem with Poncelet's theorem for quadrics
   (Section~6).

In the next section we formulate Theorem~\ref{th10} and observe its special cases for the Poncelet,
Zigzag, and Steiner theorems. In Section~3 we give a geometrical proof of Theorem~\ref{th10}.
For the sake of simplicity, in Sections 4-6 we deal with the case of {\em nested circles} $\alpha_0, \delta, \alpha_1$.
This means that $\alpha_0$ is inside $\delta$ which is inside $\alpha_1$ and all the circles
$\omega_k$ are inscribed in the annulus between $\alpha_0$ and $\alpha_1$, i.e., they are from the family~$\cM_1$.
Then, in Section~7, we prove Emch's theorem for the general position of circles. The proof remains short, but
becomes less obvious than for the nested circles.

In what follows, we denote points and vectors from $\re^2$ by bold letters,
all distances are Euclidean, the distance between points $\bx$ and $\by$
are denoted either as $\bx\by$ or as $|\bx - \by|$. By $\bc_i, r_i$, and $f_i(\bx) =
|\bx - \bc_i| - r_i^2$ we denote the center of the circle~$\alpha_i$,
its radius, and the power w.r.t. $\alpha_i$ respectively, $i=0,1$.  For two different quadrics
$\gamma_j = \{\bx \in \re^2 \ | \ q_j(\bx) = 0\},  \, j = 0,1$,  we denote
by $\{\gamma_0, \gamma_1\}$ the pencil passing through them, which is the one-parametric family of quadrics
  $\gamma_t = \{\bx \in \re^2 \ | \ (q_0 + t q_1)(\bx) = 0\},  \, t \in \bar \re$, where $\bar \re = \re\cup \{\infty\}$. For a circular series $\{\omega_k\}_{k \in \n}$, we denote by $\bx_k, \bx_{k+1}$ the points of intersection
  of the circle $\omega_k$ with $\delta$ and by $\bt_0^{k}, \bt_1^{k}$ the points of its tangency
with $\alpha_0$ and $\alpha_1$ respectively.

\medskip
\begin{center}
\textbf{2. The main result}
\end{center}
\medskip

Let $\alpha_0, \alpha_1$ and $\delta$ be arbitrary circles.
Consider a circle $\omega$ tangent to both $\alpha_0$ and $\alpha_1$ and intersecting $\delta$
at some points $\bx, \by$. Let $\omega'$ be a circle close to $\omega$ and also touching
$\alpha_0, \alpha_1$; $\bx', \by'$ be the corresponding points of intersection ($\bx'$ is close to $\bx$).
The oriented  lengths of small arcs $\stackrel{\smile }{\bx'\bx}$ and $\stackrel{\smile }{\by'\by}$ of the circle~$\delta$ as $\omega' \to \omega$ are  $d\bx$ and $d\by$. Thus, if one slightly perturbs the circle $\omega$, its points of intersection
with the circle $\delta$ moves to $d\bx$ and $d\by$.
\begin{defi}\label{d5}
Given three circles $\alpha_0, \alpha_1, \delta$ and an index~$i \in \{0,1\}$.
A Lebesgue measurable function~$\rho: \delta \to \re_+$ defines an invariant measure
if for almost all circles $\omega$ touching $\alpha_0, \alpha_1$ with index $i$ we have
\begin{equation}\label{invar}
\rho(\bx)\, |d\bx| \ = \ \rho(\by) \, |d \by|\, ,
\end{equation}
where $\bx, \by$ are points of intersection of the circles $\omega$ and $\delta$.
\end{defi}
For an arbitrary arc $\stackrel{\smile }{\bx\by} \, \subset \, \delta$ we denote by
$m(\stackrel{\smile }{\bx \by}) = \int_{\bx}^{\by} \rho(\bss)d\bss $
 its measure, or {\em mass}.
In case of nested  circles $\alpha_0, \delta, \alpha_1$, any slight perturbation of a circle $\omega$ moves the points $\bx$ and $\by$ in the same direction.
Hence,  $d\bx$ and $d\by$ always have the same sign, and equality~(\ref{invar}) becomes
$\rho(\bx)d\bx = \rho(\by)d\by$. Integrating, we obtain $m(\stackrel{\smile }{\bx\by}) \equiv {\rm const}$.
Thus, all circles $\omega \in \cM_i$
cut arcs of the same mass  $\tilde m$ from the circle $\delta$.
In particular, in Emch's theorem,  $m(\stackrel{\smile }{\bx_k\bx_{k+1}}) = \tilde m$ for all~$k$.
Hence, the circular series closes
after $n$ steps if and only if $n\, \tilde m$ is an integer multiple of $m(\delta)$.
This proves Emch's theorem in case of nested circles.
The general case is more delicate, we consider it in Section~7.
\begin{theorem}\label{th10}
The function $\rho(\bx) = \frac{1}{\sqrt{|f_0(\bx)f_1(\bx)|}}$ defines an invariant measure
on any circle~$\delta$.
\end{theorem}
Note that the function $\rho(\bx)$ is defined on the whole plane (including the circles~$\alpha_0, \alpha_1$, where it
equals to $+\infty$) and does not depend on the circle~$\delta$.
The restriction of this function to any circle defines an invariant measure on it.
Before we prove Theorem~\ref{th10} observe some of its special cases.
\smallskip

\textbf{1. The circle $\mathbf{\alpha_1}$ is infinitely big: Poncelet's theorem} If we increase the radius  of $\alpha_1$
leaving its center and all other circles unmoved, then $f_1(\bx)/r_1^2 \to 1$
uniformly on any compact subset of $\re^2$ as $r_1 \to \infty$. Hence, on
the circle $\delta$,  the function $f_1$ becomes equivalent to an identical constant.
Consequently, the function $\rho(\cdot)$ becomes proportional to $1/\sqrt{|f_0(\cdot)|}$, which is the
Jacobi-Bertrand measure. On the other hand, all the circles $\omega_k$ also enlarge
as $r_1 \to \infty$,
and their arcs touching $\alpha_0$ become close to line segments. Therefore, in the limit as $r_1 \to \infty$,
the Emch theorem becomes the
Poncelet theorem (for circles) and the measure $\rho$ becomes the Jacobi-Bertrand measure.
 Hence, the invariance property of the Jacobi-Bertrand measure
follows from Theorem~\ref{th10}.
   \smallskip

\textbf{2. The circles $\alpha_0, \alpha_1$ and $\delta$ belong to one pencil: Steiner's theorem.}
If the circle~$\delta$ belongs to the pencil $\{\alpha_0, \alpha_1\}$, then the functions $f_0$ and $f_1$ are proportional on  $\delta$:
$f_1(\bx) = - c f_0(\bx), \ \bx \in \delta$. Hence $\rho(\bx) = \frac{1}{c \, f_0 (\bx)}$.
So, in this case the reciprocal of the power w.r.t. the circle $\alpha_0$ is an
invariant measure on the circle $\delta$. If $\delta$ is the locus of points
of tangency of two circles both touching $\alpha_0$ and~$\alpha_1$, we obtain the Steiner theorem.
   \smallskip

\textbf{3. The circles $\mathbf{\alpha_0}$ and  $\mathbf{\alpha_1}$ are concentric: the Zigzag theorem.}
If $\alpha_0$ and $\alpha_1$ are concentric, then the Emch theorem becomes the Zigzag theorem for the circles $\delta$ and $\alpha$ (the circle $\alpha$ is of radius $r = \frac{r_0+r_1}{2}$ and is concentric to
$\alpha_0, \alpha_1$) and for the jump length
$l = \frac{|r_1-r_0|}{2}$. The measure $\rho(\cdot)$ on the circle $\delta$ becomes the
Black-Howland measure $b(\cdot)$ for Zigzag theorem~\cite{BHH}. It is defined as $b(\bx) = 1/|(\bx - {\bc}_0)\times ({\bx - \bz})|$,
where $\times$ denotes the operation of cross (vector) product, $\bx \in \delta$ and $\bz \in \alpha$ is such that $|\bx - \bz| = l$. In other terms, $1/b(\bx)$ is the double area
of a triangle with the sidelengths $x = |\bx - \bc_0|, \, \frac{r_1 + r_0}{2}$ and $\frac{|r_1 - r_0|}{2}$
(fig. 2).

\begin{figure}[htb]
\center
\includegraphics[scale=0.45]{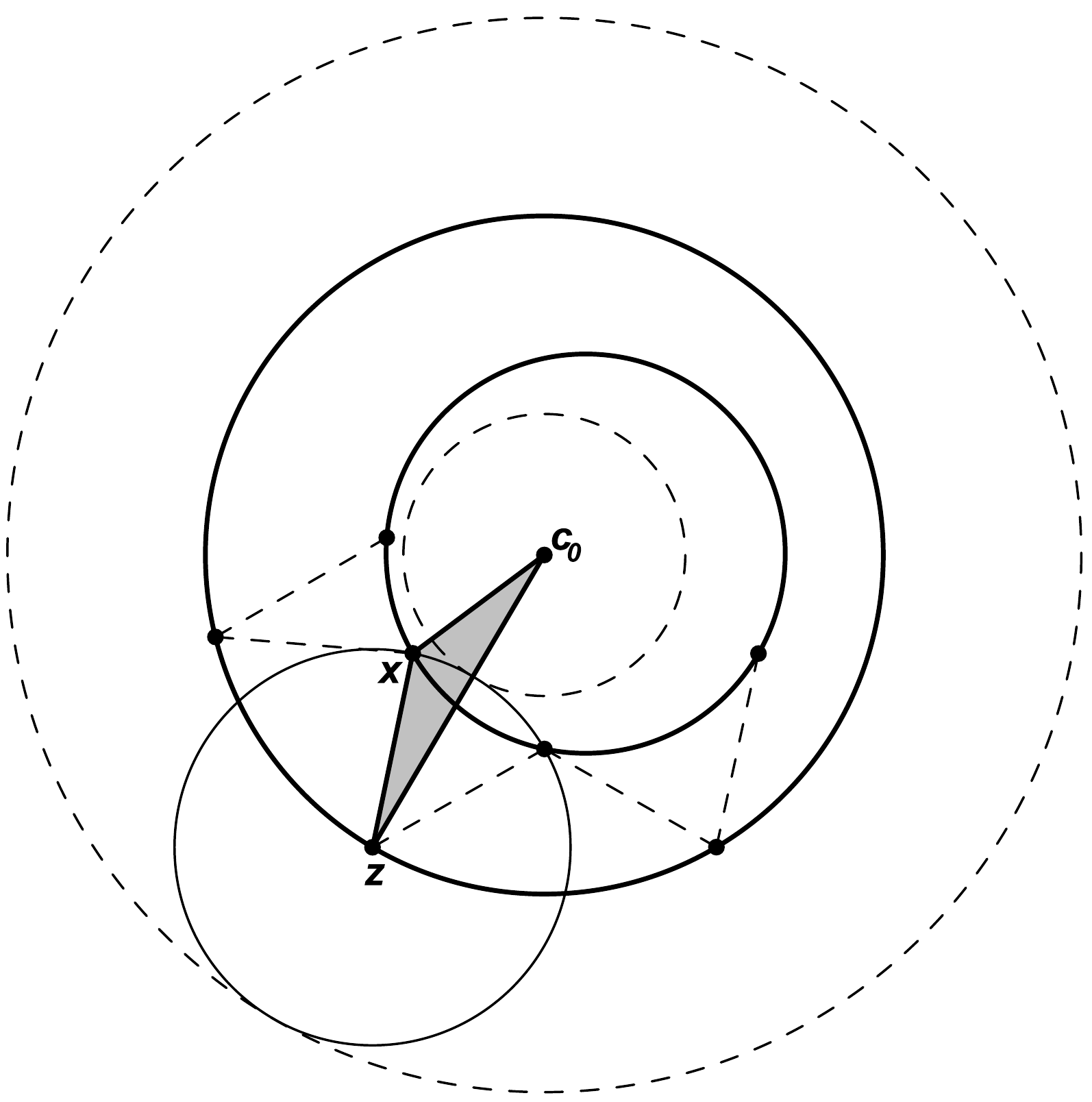}
\caption{{\footnotesize The invariant measure for the Zigzag theorem}}
\label{f2}
\end{figure}

The Heron formula yields
$$
\frac{1}{b(\bx)} = \frac12 \sqrt{(r_1+x)(r_1-x)(x+r_0)(x-r_0)} = \frac{1}{2}\sqrt{(r_1^2 - x^2)(x^2 - r_0^2)} =
\frac{1}{2}\sqrt{- f_1 (\bx) \cdot f_0 (\bx)}\, .
$$
 Hence $b(\bx) = 2\rho(\bx)$ for all $\bx \in \delta$. Thus, the Black-Howland measure is the special case of $\rho(\cdot)$, when the circles
$\alpha_0$ and $\alpha_1$ are concentric.

\bigskip

\begin{center}
\textbf{3. Proof of the main theorem}
\end{center}
\bigskip

Let a circle $\omega$ touch both $\alpha_0$ and $\alpha_1$. We consider the line
connecting the two tangent points and denote by
$h(\bx)$ the distance from a point $\bx$ to that line. We are going to show that
the function $\rho$ on the circle $\omega$ is proportional to $1/h$.
\begin{prop}\label{p10}
Suppose $\omega$ is an arbitrary circle touching $\alpha_0$ and $\alpha_1$
at points ${\bt}_0$ and ${\bt}_1$ respectively; then the restriction of the function
$\rho(\bx) = {1}/{\sqrt{|f_0(\bx)f_1(\bx)|}}$ to $\, \omega$ is proportional to the reciprocal of the distance to the line~${\bt}_0{\bt}_1$. Thus,
$\, \rho(\bx) \, \sim \, 1/h(\bx), \ \bx \in \omega$.
\end{prop}
{\tt Proof.} Let the line $\bx{\bt}_0$ meet the circle $\alpha_0$ for the second time at point ${\bz}_0$.
Note that $f_0(\bx) = \bx {\bt}_0\cdot \bx {\bz}_0 = c\cdot \bx {\bt}_0^2$, where $c$ is a constant. Indeed, since
the circles $\alpha_0$ and $\omega$ are homothetic with respect to the point of tangency~${\bt}_0$,
 the ratio ${\bz}_0{\bt}_0/ \bx {\bt}_0$ is constant, and hence so is the ratio $\bx {\bz}_0/ \bx {\bt}_0$.
Similarly, $f_1(\bx)$ is proportional to $(\bx {\bt}_1)^2$. Thus, $\sqrt{|f_0(\bx)\cdot f_1(\bx)|}
\sim \bx {\bt}_0 \cdot \bx {\bt}_1$, which is proportional to the area of the triangle~$\triangle\, {\bt}_0\bx {\bt}_1$
(because $\sin (\angle \, {\bt}_0\bx {\bt}_1)$ is constant), which is, in turn, proportional to
its altitude $h(\bx)$, since this triangle has a constant base~${\bt}_0{\bt}_1$.

   {\hfill $\Box$}
\medskip

A different proof of Proposition~\ref{p10} based on properties of
pencils of quadrics is given in Section~5, where we prove  a generalization
of Emch's theorem.

\begin{prop}\label{p20}
Suppose a circle $\omega$ passes through points $\bk$ and  $\bl$ and
meets a circle $\delta$ at points $\bx$ and $\by$; then a small
perturbation of $\omega$ that passes through $\bk$ and $\bl$ satisfies
$\frac{|d\by|}{|d\bx|} = \frac{q(\by)}{q(\bx)}$, where $q(\cdot)$ is the distance to the line $\bk \bl$.
\end{prop}
{\tt Proof.} Since three pairwise chords of three circles concur,
the lines $\bx \by$ and $\bx' \by'$ meet on the line $\bk \bl$, at some point $\bn$.
Equalities $\bx' \bn\cdot \by' \bn = \bk \bn \cdot \bl \bn = \bx \bn \cdot \by \bn$ imply similarity of triangles
$\triangle\, \bx \bn \bx' \, \sim \, \triangle\, \by' \bn \by$, which yields $\frac{\by \by'}{\bx \bx'} = \frac{\bn \by'}{\bn \bx}$.
Replacing $\bn \bx'$ by a close value $\bn \bx$ and $\frac{\bn \by}{\bn \bx}$ by $\frac{q(\by)}{q(\bx)}$, we conclude the proof.

   {\hfill $\Box$}
\medskip

\begin{figure}[htb]
\center
\includegraphics[width=0.45\textwidth]{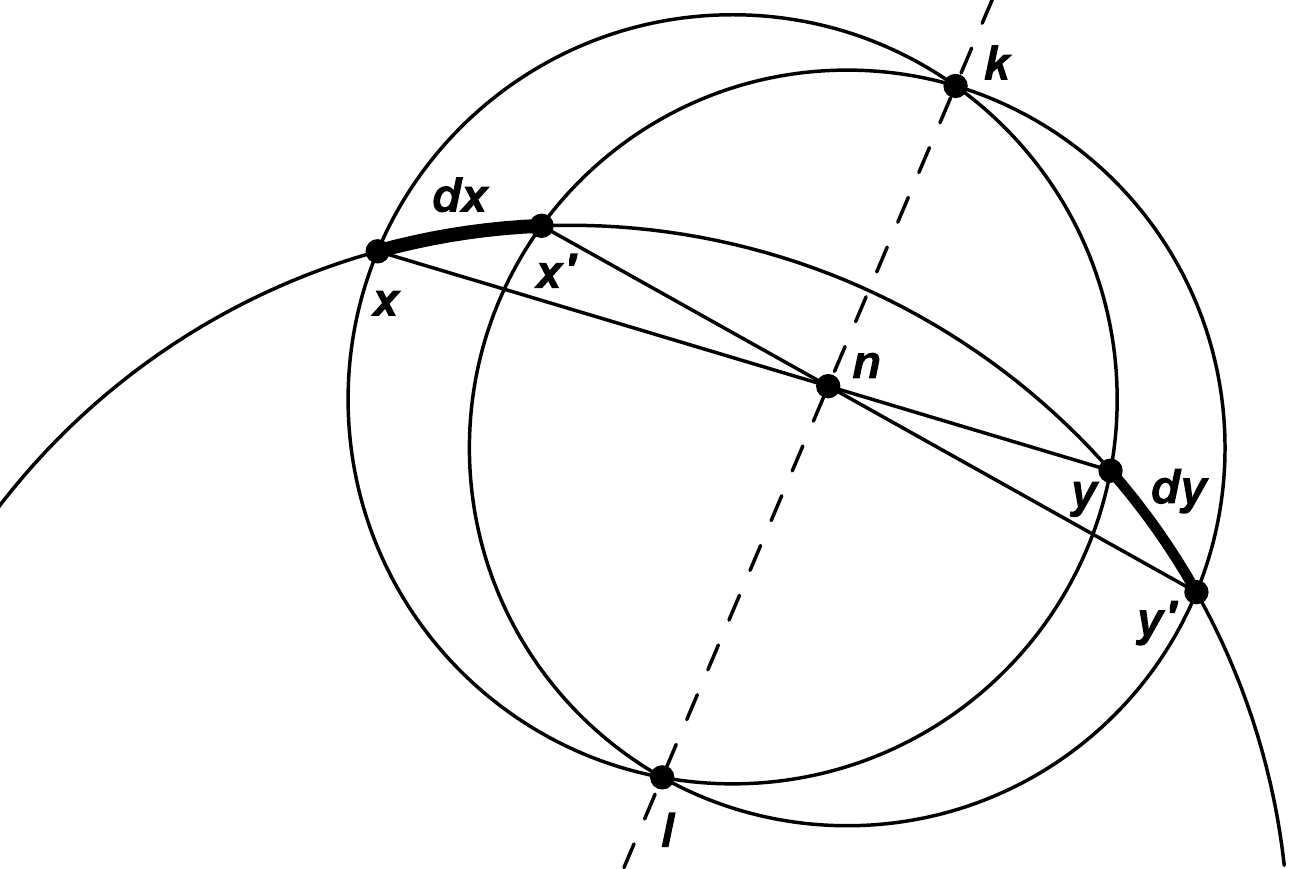}
\caption{{\footnotesize Proof of Proposition 2}}
\label{f3}
\end{figure}

{\tt Proof of Theorem~\ref{th10}.} Let $\omega$ be a circle
touching  $\alpha_0, \alpha_1$ and intersecting $\delta$ at points $\bx, \by$;
$\omega'$ be a small perturbation of $\omega$. Replacing in the equality~(\ref{invar}) the function $\rho(\cdot)$
by $1/h(\cdot)$ (Proposition~\ref{p10}),
$h(\cdot)$ by a close value $q(\cdot)$, which is the distance to the common chord of the circles $\omega$ and $\omega'$
(this chord tends to the line ${\bt}_0{\bt}_1$, hence $q/h \to 1$ as $\omega' \to \omega$) we come to
an equivalent assertion $\frac{|d \bx |}{q(\bx)} =  \frac{|d \by|}{q(\by)}$,
 which follows from Proposition~\ref{p20}.

   {\hfill $\Box$}
\medskip

\begin{center}
\textbf{4. Generalizations to pencils of circles }
\end{center}
\bigskip

For the sake of simplicity, in Sections 4-6 we consider the case of nested circles $\alpha_0, \delta, \alpha_1$.
The general case is analysed in Section~7.

The measure $\rho$ provides a simple way to generalize Emch's theorem from one pair of circles
$(\alpha_0, \alpha_1)$ to arbitrary sequence of pairs $(\alpha_0^{(k)}, \alpha_1^{(k)})_{k \in \n}$,
where each $\alpha_i^{(k)}$ is taken from a given pencil of circles~$\cA_i$.
Such an extension for Poncelet's theorem is known, it was  proved by Poncelet himself~\cite{Pons}, then
developed by Lebesgue~\cite{L}, see also~\cite{B}. A similar extension for Emch's
theorem originated in~\cite{P2}. Let $\cA_0, \cA_1$ be arbitrary pencils of circles both containing
the circle~$\delta$. Take arbitrary sequences $\{\alpha_0^{(k)}\}_{k \in \n}\subset \cA_0$ and
$\{\alpha_1^{(k)}\}_{k \in \n}\subset \cA_1$.
\begin{prop}\label{p30}
  All the pairs $(\alpha_0^{(k)}, \alpha_1^{(k)}), \ k \in \n$, generate invariant measures on the circle~$\delta$ that
  are proportional to one measure $\rho$.
\end{prop}
{\tt Proof.} Let $f^{(k)}_i$
denote the power w.r.t. the circle $\alpha^{(k)}_i, \ i = 0,1$.
Since this circle belongs to the pencil
$\{\delta, \alpha^{(1)}_i\}$, it follows that $f^{(k)}_i =
(1-t_{i,k}) f_{\delta} + t_{i,k} f^{(1)}_i$, for some $t_{i,k} \in \bar \re$.
For all $\bx \in \delta$, we have $f_{\delta}(\bx) = 0$, and hence
$f^{(k)}_0(\bx)f^{(k)}_1(\bx) = t_{0, k}t_{1,k} f^{(1)}_0(\bx)f^{(1)}_1(\bx)$,
i.e., the measures generated by the $k$th pair and by the first pair are proportional on~$\delta$.

   {\hfill $\Box$}
\medskip

Thus, for given pencils~$\cA_0, \cA_1$ containing a circle $\delta$,
every pair $(\alpha_0, \alpha_1) \in \cA_0\times \cA_1$ generates  an invariant  measure, and all
those measures are proportional on $\delta$.
Hence, the following generalized Emch's theorem holds.
Let $(\alpha_0^{(k)}, \alpha_1^{(k)}) \in \cA_0\times \cA_1, \ k \in \n$,  be an arbitrary sequence of pairs.
Consider a circular series $\{\omega_k\}_{k \in \n}$,
 where $\omega_k$ touches the $k$th pair (fig. 4).
\begin{figure}[htb]
\center
\includegraphics[scale=0.45]{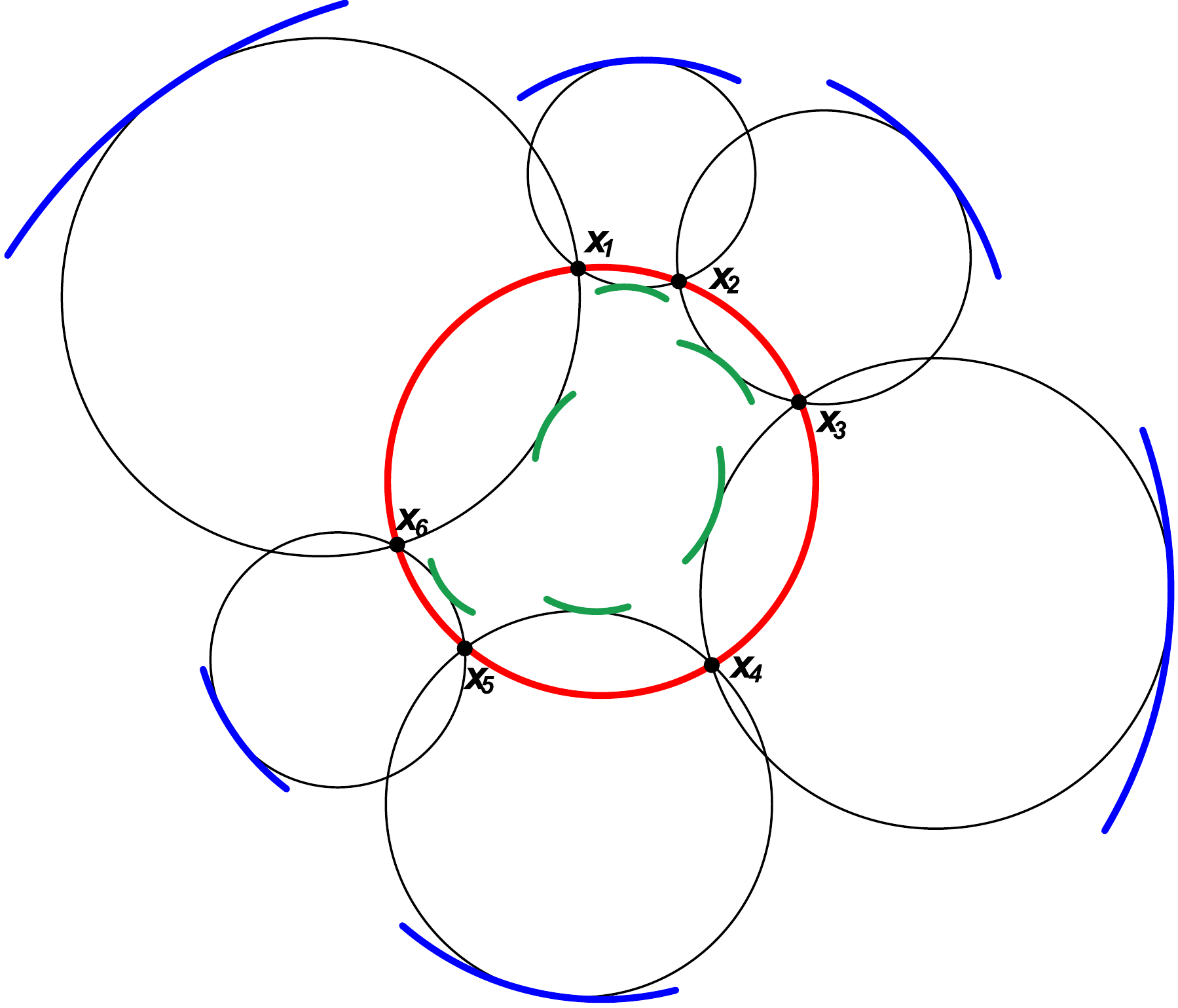}
\caption{{\footnotesize Emch theorem for pencils of circles}}
\label{f4}
\end{figure}

 If for some initial circle $\omega_1$, we have
$\omega_{n+1} = \omega_1$, then it holds for arbitrary $\omega_1$ touching the first pair. Moreover,
 after an arbitrary change of order of the pairs~$(\alpha_0^{(1)}, \alpha_1^{(1)}), \ldots ,
 (\alpha_0^{(n)}, \alpha_1^{(n)})$, this series still closes after $n$ steps. See~\cite{P2}
 for the precise formulation. The proof is literally the same as for the Emch theorem.
  The closing after $n$ steps takes place if and only if the sum of masses of $n$ arcs
  $\stackrel{\smile }{\bx_1\bx_2}, \ldots , \stackrel{\smile }{\bx_n\bx_{n+1}}$ of the circle $\delta$
 cut by the circles $\omega_k$ is equal to  a multiple of the total mass of~$\delta$. This equality depends neither on
 the location of the initial circle $\omega_1$ (due to the invariance of the measure) nor on the
 ordering of the circles (due to commutativity of summation).

 Several corollaries can be drawn from Proposition~\ref{p30} even if
 the circular series does not close.
  They are based on the following simple observation.
  \begin{prop}\label{p40}
Under the assumptions of Emch's theorem, for every $\tilde m > 0$,
the following holds: all circles $\omega$ that cut from $\delta$ arcs of the same mass $\tilde m$
(generated by the measure $\rho(\bx) = 1/\sqrt{|f_0(\bx)f_1(\bx)|}$)
and touch  $\alpha_0$ with a given index,
touch a fixed circle from the pencil $\cA_1 = \{\delta, \alpha_1\}$.
\end{prop}
{\tt Proof.} For an arbitrary circle $\omega$, the pencil $\cA_1$ contains a unique circle $\alpha_1'$ touching
$\omega$ with a given index. By Proposition~\ref{p30}, the measure $\rho$
is invariant for the pair $(\alpha_0, \alpha_1')$, hence all
circles $\omega'$ touching this pair with a given index cut the same mass $m$~on $\delta$.

   {\hfill $\Box$}
\medskip

\begin{cor}\label{c10}
Let us have two circular series $\{\omega_k\}_{k \in \n}$ and $\{\omega_k'\}_{k \in \n}$ touching circles $\alpha_0, \alpha_1$ with the same index and having the same direction. Let $\omega_k$ and $\omega_k'$ intersect the circle $\delta$ at points $\bx_k, \bx_{k+1}$ and $\bx_k', \bx_{k+1}'$ respectively. Denote by $\gamma_k$ the circle passing through the points
$\bx_k$ and $\bx_k'$ and touching $\alpha_0$ with the same index. Then
all $\gamma_k$ touch a fixed circle from the pencil $\cA_1 = \{\delta, \alpha_1\}$.
\end{cor}
{\tt Proof.} All the arcs $\stackrel{\smile }{\bx_k\bx_k'}$ of the circle $\delta$ have the same masses. Invoking Proposition~\ref{p40}
completes the proof.

   {\hfill $\Box$}
\medskip

\begin{cor}\label{c20}
Let $\{\omega_k\}_{k \in \n}$ be a circular series touching $\alpha_0, \alpha_1$
and let $\omega_k$ intersect the circle $\delta$ at points $\bx_k, \bx_{k+1}$.
Fix $r \in \n$ and for every $k$, consider a  circle passing through
$\bx_k$ and $\bx_{k+r}$  and touching $\alpha_0$ with the same index. Then
all those circles touch a fixed circle $\alpha_r \in \cA_1$.
\end{cor}
{\tt Proof.} We apply Corollary~\ref{c10} with $\omega_k' = \omega_{k+r}$.

   {\hfill $\Box$}
\medskip

Thus, the situation is the same as for the diagonals of Poncelet's polygons~\cite{B}. Here, if
a curvilinear broken line is inscribed in a circle $\delta$ and its sides touch a pair of circles
$\alpha_0, \alpha_1$, then all its diagonals of $r$th order touching $\alpha_0$ also touch a fixed circle
from the pencil~$\cA_1 = \{\delta, \alpha_1\}$.

\bigskip

\begin{center}
\textbf{5. Emch's theorem for cyclics}
\end{center}
\bigskip

A {\em cyclic} is a plane algebraic curve of order four
defined by the equation
\begin{equation}\label{cyc}
F(x_1, x_2) \ = \ \lambda (x_1^2 + x_2^2)^2 \, + \, (x_1^2 + x_2^2) \, \ell (x_1, x_2) \, + \, Q(x_1,x_2) \ = \ 0\, ,
\end{equation}
where $\ell$ is a linear form and $Q$ is a polynomial of degree at most two.
In the sequel we assume $\lambda = 1$;
the general case follows  either by normalization (if $\lambda \ne 0$) or
by a limit passage (if $\lambda = 0$).
A pair of circles on the plane is always a cyclic, but not vice versa. An arbitrary quadric is
a cyclic as well.
Some properties of cyclics can be found in~\cite[chapter 4, sect. 2]{C}.
Nilov in~\cite{N} proved that the Emch theorem remains true after replacing  the pair of circles
$\alpha_0, \alpha_1$ by an arbitrary cyclic~$\Gamma$. In this case, all circles $\omega_k$ have
{\em double tangency} (i.e. two points of tangency) with~$\Gamma$.

\begin{figure}[htb]
\center
\includegraphics[scale=0.4]{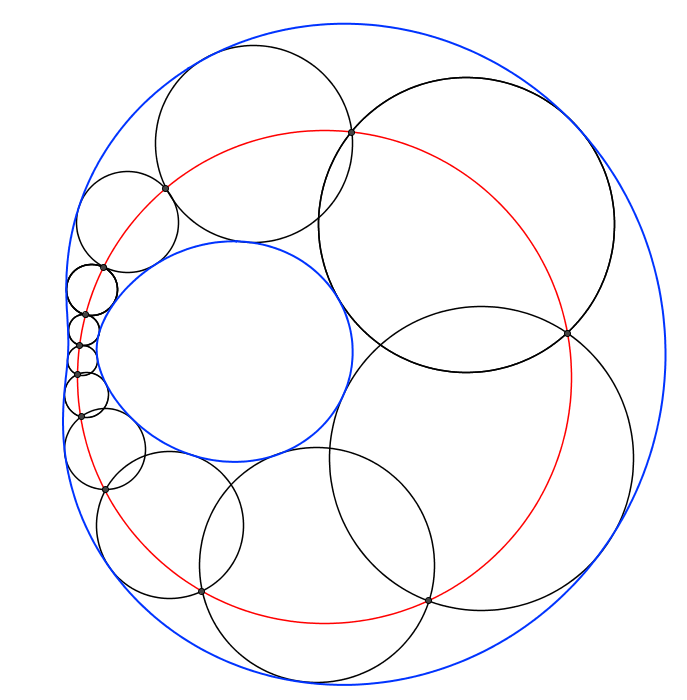}
\caption{{\footnotesize The Emch theorem for a cyclic}}
\label{f5}
\end{figure}

If we take into account the complex tangency,
then there exist four families
of circles with double tangency with~$\Gamma$, all $\omega_k$ belong to one of them~\cite{C}.
The definition of invariant measure remains the same. The proof in~\cite{N} is geometrical and relies
 on the Poncelet  theorem for quadrics. The universal measure $\rho$ enables us to
 give a self-contained  proof using the following generalization of Theorem~\ref{th10}.
  \smallskip

 \noindent \textbf{Theorem~\ref{th10}'}. {\em The function $\rho(\bx) = 1/\sqrt{|F(\bx)|}$ generated by a cyclic
 $ {\Gamma = \{\bx \in \re^2  \, | \, F(\bx) = 0\}}$  defines an invariant measure on any circle.}
 \smallskip

The proof is based on the following generalization of
Proposition~\ref{p10} (Section~3) to cyclics.

\smallskip

 \noindent \textbf{Proposition~\ref{p10}'}. {\em If a circle $\omega$ touches a given cyclic~$\Gamma = \{\bx \in \re^2  \, | \, F(\bx) = 0\}$ at points
 ${\bt}_0, {\bt}_1$, then the function
$F|_{\omega}$  is proportional to the square of the distance to the line~${\bt}_0{\bt}_1$. }
 \smallskip

Combining this with Proposition~\ref{p20} we prove Theorem~\ref{th10}' in the same way as Theorem~\ref{th10}.
The proof of Proposition~\ref{p10}' uses an algebraic argument and one auxiliary result.
\begin{defi}\label{d10}
For an arbitrary circle $\delta$, two algebraic curves are called $\delta$-equivalent, if the polynomials defining those curves are proportional on~$\delta$.
\end{defi}
The $\delta$-equivalence of curves $g_1(\bx) = 0$ and $g_2(\bx) = 0$
means that for some $\mu \ne 0$, the polynomial $g_1 - \mu g_2$ is divisible by $f_{\, \delta}$
(the power w.r.t. $ \delta$).
\begin{lemma}\label{l10}
For an arbitrary cyclic~$\Gamma$ and a circle~$\delta$, the closure of the set of quadrics
$\delta$-equivalent to~$\Gamma$ is a pencil of quadrics containing~$\delta$.
\end{lemma}
{\tt Proof.} Let
$f_{\, \delta} (\bx) = x_1^2 + x_2^2 + \ell_{\delta} (x_1, x_2) + A_{\delta} = 0$, where
$\ell_{\delta}$ is a linear form and $A_{\delta}$ is a constant. A polynomial $p(\bx)$
possesses the property ${\rm deg} \, (F \, - \, p \, f_{\, \delta}) \, \le \, 2$ if and only if
\begin{equation}\label{q-cond}
p (\bx) \   = \
x_1^2 + x_2^2 + \ell_p (x_1, x_2) + A_{p} \, , \qquad \mbox{with} \qquad \ell_p + \ell_{\delta} = \ell\, ,
\quad A_{p} \in \re\, ,
\end{equation}
where $\ell$ is from equation of cyclic~(\ref{cyc}).
Denote by $p_0$ the
polynomial~(\ref{q-cond}) with $A_p = 0$. For arbitrary $A_p \in \re$,  we have a quadratic polynomial
$F \, - \, p \,  f_{\, \delta} \, = \, F \, - \, (p_0 + A_p) f_{\, \delta} \, = \,
\bigl( F \, - \, p_0 \,  f_{\, \delta}\bigr) \, - \, A_p f_{\, \delta}$.
When $A_p$ runs over  $\bar \re $, these polynomials define a pencil of quadrics
which contains $\delta$ (for $A_p = \infty$).

   {\hfill $\Box$}
\medskip

{\tt Proof of Proposition~\ref{p10}'.} Any quadric $\omega$-equivalent to $\Gamma$
touches the circle $\omega$ at points $\bt_0$ and $\bt_1$. By Lemma~\ref{l10},
those quadrics form a pencil.
On the other hand, all quadrics touching a circle at two  points form a pencil
that contains a double line connecting those points~\cite[section 16.4.10]{B}. Hence, these two pencils
coincide. In particular,
the double line $\bt_0\bt_1$ is $\omega$-equivalent to $\Gamma$. So, the function $F|_{\, \omega}$
 is proportional to the square of the distance to the line $\bt_0\bt_1$.

   {\hfill $\Box$}

\bigskip

\newpage
\begin{center}
\textbf{6. The Emch theorem  and Poncelet's theorem for quadrics}
\end{center}
\bigskip

By Lemma~\ref{l10}, a cyclic is equivalent to a quadric on every circle. Moreover,
if a cyclic~$\Gamma$ and a circle~$\delta$ are fixed, then all such quadrics form a pencil~$\cQ$.
This implies that the invariant measure~$\rho = 1/\sqrt{|F|}$ generated by~$\Gamma$
on the circle $\delta$ coincides with the Jacobi-Bertrand measure $1/{\sqrt{|q|}}$ generated by any
quadric from~$\cQ$. Therefore, $\cQ$ contains a quadric~$\gamma$ tangent to all
 lines $\bx_{k}\bx_{k+1}, \, k \in \n$, corresponding to a circular series~$\{\omega_k\}$.
 Hence the   Emch  theorem follows from Poncelet's theorem for quadrics $\delta$ and $\gamma$.

 It was first noted by Hrask\'o~\cite{H} that the Zigzag theorem  can be derived from
 Poncelet's theorem for quadrics. Then in~\cite{P1} this result was extended to Emch's theorem, and
 in~\cite{N} to cyclics.
 The proofs in those works are different and nontrivial.
Now we see that this is actually  a consequence of
equivalence of a cyclic to a certain quadric on a circle. Moreover, it is possible to find the
 desired quadric~$\gamma$ explicitly. We have $F(\bx) = f_0(\bx)f_1(\bx)$, where
$f_i (\bx) = x_1^2 + x_2^2 + \ell_i(x_1, x_2) + B_i = 0$ is the power w.r.t. the circle~$\alpha_i, \, i = 0,1$. Applying~(\ref{q-cond})
we see that the polynomial $p(\bx) = x_1^2 + x_2^2 + \ell_p (x_1, x_2) + A_p$
satisfies the equalities $\, \ell_p + \ell_{\, \delta} \, =\, \ell_0 + \ell_1 $.
The quadric $\gamma$ is thus given by the equation $q(\bx)  =
(f_0f_1 -  p f_{\, \delta})(\bx) = 0$. Simplifying, we get
\begin{equation}\label{q}
q(\bx) \ = \ \bigl(\ell_0(\bx) + B_0\bigr)\, \bigl(\ell_1(\bx) + B_1\bigr) \ -\
\bigl(\ell_{\, \delta}(\bx) + A_{\delta}\bigr)\, \bigl(\ell_p (\bx) + A_p\bigr)\ + \
\bigl(x_1^2 + x_2^2\bigr)\, \bigl(B_0+B_1 - A_{\delta} - A_p\bigr),
\end{equation}
where $\ell_p = \ell_0 + \ell_1 - \ell_{\delta}$,
and the parameter $A_p$ is found by the tangency condition.

The inverse implication can also be easily realized. If we have a circle~$\delta$ and a quadric~$\gamma$,
then one can find functionals $\ell_0, \ell_1, \ell_p$ and constants $B_0, B_1, A_p$
such that $\ell_0 + \ell_1 = \ell_p + \ell_{\delta}$ and~(\ref{q}) holds.
This way  we find circles $\alpha_0, \alpha_1$
such that all chords $\bx_{k}\bx_{k+1}, \, k \in \n$, in the Emch theorem are tangent to $\gamma$.
Hence, Emch's theorem implies the Poncelet's theorem for a circle and a quadric, which is
equivalent to the case of two quadrics (by means of a suitable stereographic
projection).

Thus,  {\em the Poncelet theorem for quadrics follows from
Emch's theorem}.

\bigskip

\begin{center}
\textbf{7. The Emch theorem for general position of circles}
\end{center}
\bigskip

As we noted in Section~2, the very existence of an invariant measure immediately implies the Emch theorem
for nested circles. In this case, the differentials $d\bx$ and $d\by$ in~(\ref{invar})
always have the same sign. In particular, for a small perturbation of the circular series
$\{\omega_k\}$, we have ${\rho(\bx_k) d\bx_k \equiv {\rm const}, \, k \in \n}$.
Integrating, we obtain that if the circle $\omega_1$ moves to a circle $\omega_1'$,
then for the series $\{\omega_k\}$ and $\{\omega_k'\}$, we have
$m(\stackrel{\smile }{\bx_k\bx_{k}'}) \equiv {\rm const}, \, k \in \n$.
In particular, $m(\stackrel{\smile }{\bx_1\bx_{1}'}) = m(\stackrel{\smile }{\bx_{n+1}\bx_{n+1}'})$, hence if $\bx_{n+1} = \bx_1$, then
 $\bx_{n+1}' = \bx_1'$, which completes the proof. In the general case, however,
a small perturbation of $\omega_1$
can move the points $\{\bx_k\}_{k \in \n}$ in different directions.
 That is why, to prove Emch's theorem in the general case
 we need to modify the invariant $\rho(\bx)|d\bx|$ to respect  the sign of the differential~$d\bx$.

For an arbitrary triangle $\triangle \, \ba \bb \bc$ we denote by $\tau(\ba \bb \bc)$
its orientation: $\tau(\ba \bb \bc) = 1$ if its vertices follow in the positive direction or,
equivalently, the pair of vectors $\bb - \ba$ and $\bc - \ba$ is positively oriented. Otherwise,
$\tau(\ba \bb \bc) = - 1$. To avoid considering two cases, we make the following assumption:
\smallskip

\noindent \textbf{Assumption 1.} {\em The circle $\omega_1$ lies inside $\alpha_1$.}

\smallskip

This assumption is not restrictive, it can always be achieved by a suitable
inversion.  Note also that if $\omega_1$ lies inside $\alpha_1$, then so does $\omega_2$
(since it intersects $\omega_1$), and $\omega_3$, etc. Thus, Assumption~A means that the whole
series $\{\omega_k\}$ is inside $\alpha_1$.

\begin{theorem}\label{th20}
For any circles $\alpha_0, \alpha_1, \delta$ and for an arbitrary
circular series $\{\omega_k\}_{k \in \n}$ touching $\alpha_0, \alpha_1$ and satisfying Assumption 1, we have
$\ \tau (\bx_k \bt_0^k \bt_1^k)\, \rho(\bx_k) \, d\bx_k \, \equiv \, {\rm const}, \  k \in \n$.
\end{theorem}
    The proof of Theorem~\ref{th20} requires two auxiliary facts. The first one is a generalization
of Proposition~\ref{p20}.
\begin{prop}\label{p50}
Under the assumptions of Proposition~\ref{p20}, we have
$\frac{d\by}{d\bx} \, = \, -\, \frac{\tau(\by\bk\bl)}{\tau(\bx\bk\bl)}\, \frac{q(\by)}{q(\bx)}$.
\end{prop}
{\tt Proof.} If the chords $\bk \bl$ and $\bx \by$ intersect, then
the arcs $\stackrel{\smile }{\bx \bx'}$ and $\stackrel{\smile }{\by \by'}$ have the same sign and
$\frac{\tau(\bx\bk\bl)}{\tau(\by\bk\bl)} = -1$. Otherwise
those arcs have opposite signs and $\frac{\tau(\bx\bk\bl)}{\tau(\by\bk\bl)} = -1$.
Applying Proposition~\ref{p20}, we conclude the proof.

   {\hfill $\Box$}
\medskip

The proof of the following fact is elementary and we omit it.
\begin{lemma}\label{l20}
Let circles $\omega$ and $\nu$ pass through a point $\bm$,
circles $\alpha_0$ and $\alpha_1$ touch them with index~$0$
at points $\bt_0, \bt_1$ and $\bss_0, \bss_1$ respectively. Then
$\tau(\bm \bt_0 \bt_1) = - \tau(\bm \bss_0 \bss_1)$.
\end{lemma}
{\tt Proof of Theorem~\ref{th20}.} Arguing as in the proof of Theorem~\ref{th10} and using
Proposition~\ref{p50} for $\bx = \bx_{k}, \by = \bx_{k+1}$,
we obtain $\tau (\bx_k \bt_0^k \bt_1^k)\rho(\bx_k)
d\bx_k = - \tau (\bx_{k+1} \bt_0^k \bt_1^k)\rho(\bx_{k+1})
d\bx_{k+1}$. Applying now Lemma~\ref{l20} to the circles
$\omega = \omega_k, \, \nu = \omega_{k+1}$ and taking into account that
$\alpha_0$ and $\alpha_1$ touch them with index~$0$, because $\omega_k$ lies inside
$\alpha_1$, we conclude $\tau (\bx_{k+1} \bt_0^k \bt_1^k) = - \tau (\bx_{k+1} \bt_0^{k+1} \bt_1^{k+1})$.

   {\hfill $\Box$}

\noindent Now we are ready to prove Emch's theorem in general case.
\smallskip

{\tt Proof of the Emch theorem.} Consider a perturbation of
the circular series $\{\omega_k\}$ that moves it to a series $\{\omega_k'\}$.
The orientation of all triangles
$\triangle \bx_{k}\bt_0^k\bt_1^k, \, k = 1, \ldots , n+1$,  is not changed, whenever the perturbation
is small enough. If $\omega_{n+1} = \omega_1$,
    then the points $\bx_{n+1},  \bt_0^{n+1},  \bt_1^{n+1}$ coincide with
    $\bx_{1},  \bt_0^{1},  \bt_1^{1}$ respectively. Hence
    $\tau (\bx_{n+1} \bt_0^{n+1}  \bt_1^{n+1}) = \tau (\bx_{1} \bt_0^{1}  \bt_1^{1})$,
    and therefore  $\rho(\bx_{n+1})d\bx_{n+1} = \rho(\bx_1)d\bx_1$.
    Integrating, we obtain $m(\bx_{n+1}\bx_{n+1}') = m(\bx_{1}\bx_1')$, hence $\bx_{n+1}' = \bx_1'$.
We see that the assertion $\bx_{n+1}' = \bx_1'$ is locally stable (under small perturbations).
The continuity implies that it holds identically.

   {\hfill $\Box$}


\medskip

\begin{thebibliography}{}

\bibitem{A}
E.A.Avksentyev,
\newblock {\em A universal measure for a pencil of conics and the Great Poncelet theorem}, \,
\newblock  Sb. Math., 205 (2014), no 5, 613--632.
\smallskip

\bibitem{A2}
E.A.Avksentyev,
\newblock {\em The Great Emch Closure Theorem and a combinatorial proof of Poncelet's Theorem}, \,
\newblock  Sb. Math., 206 (2015), no 11, 1509–-1523.
\smallskip

\bibitem{BB}
W.Barth and Th.Bauer,
\newblock {\em Poncelet theorems},
\newblock  Expositiones Mathematicae, 14 (1996), 125--144.
\smallskip

\bibitem{B}
M.Berger,
\newblock {\em G\'eom\'etrie},
\newblock CEDIC, Paris (1977).
\smallskip

\bibitem{BHH}
W.L.Black, H.C.Howland and B.Howland,
\newblock {\em A theorem about zigzags between two circles},
\newblock  Amer. Math. Monthly, 81 (1974), 754--757.
\smallskip

\bibitem{Bot}
O.Bottema,
\newblock {\em Ein Schliessungssatz f\"ur zwei Kreise},
\newblock  Elem. Math., 20 (1965), 1--7.
\smallskip

\bibitem{C}
J.L.Coolidge,
\newblock {\em A treatise on the circle and the sphere, by Julian Lowell Coolidge},
\newblock Oxford: Clarendon Press, 1916.
\smallskip

\bibitem{DR}
V.Dragovi\'c and M.Radnovi\'c,
\newblock {\em Poncelet porisms and beyond},
\newblock Frontiers in Math., Birkhauser/Springer Basel AG, Basel, 2011.
\smallskip

\bibitem{E}
A.Emch,
\newblock {\em An application of elliptic functions to
Peaucellier's link-work (inversor)},
\newblock Ann.Math., ser. 2, vol. 2 (1901), 60--63.
\smallskip

\bibitem{F}
L.Flatto,
\newblock {\em Poncelet's theorem},
\newblock AMS, Providence, RI, 2009.
\smallskip

\bibitem{HH}
L.Halbeisen and N.Hungerb\"uhler,
\newblock {\em A simple proof of Poncelet's theorem
(on the occasion of its bicentennial)},
\newblock Amer. Math. Monthly, 121 (2014), no 1, 1--14.
\smallskip

\bibitem{H}
 A.Hrask\'o,
\newblock {\em Poncelet-type problems, an elementary approach},
\newblock Elem.Math., 55 (2000), 45--62.
\smallskip

\bibitem{K}
 V.V.Kozlov,
 \newblock {\em Rationality conditions for the ratio of elliptic integrals and the great Poncelet theorem},
 \newblock  Moscow Univ. Math. Bull., 58 (2003), no 4, 1--7.
\smallskip

\bibitem{L}
H.Lebesque,
\newblock {\em Les Coniques},
\newblock  Cauthier-Villars, Paris (1942).
\smallskip

\bibitem{N}
F.Nilov,
\newblock {\em Families of conics and circles with double tangencies},
\newblock  Sb. Math., submitted.
\smallskip

\bibitem{Pons}
J.V.Poncelet,
\newblock {\em Trait\'e des propri\'et\'es projectives des figures},
\newblock Paris 1865, (first ed. in 1822).
\smallskip

\bibitem{P1}
V.Yu.Protasov,
\newblock {\em One generalization of Poncelet's theorem},
\newblock Russian Math. Surveys, 61 (2006), no 6, 187--188.
\smallskip

\bibitem{P2}
V.Yu.Protasov,
\newblock {\em Generalized closing theorems},
\newblock Elem. Math.,  (2011), 66 (2011), no 3, 98--117.
\smallskip

\bibitem{Sc}
I.J.~Schoenberg,
\newblock {\em On Jacobi-Bertrand's Proof of a Theorem of Poncelet},
\newblock Studies in Pure Mathematics, Birkhauser, Boston, 1983, 623--627.
\smallskip

\end{thebibliography}
\end{document}